\documentclass[letterpaper]{amsart} 

\usepackage{amsmath}
\usepackage[all]{xy}

\usepackage{psfrag}
\usepackage{epsfig}
 
\numberwithin{equation}{section}

\newtheorem{theorem}{Theorem}[section] 
\newtheorem{proposition}[theorem]{Proposition}

\newtheorem{lemma}[theorem]{Lemma} 
 
\theoremstyle{definition}

\newtheorem{definition}[theorem]{Definition}

\def\proof{\smallskip\noindent {\bf Proof.\ }} 
\def\endproof{\hfill$\square$\medskip} 

\def\ZZ{\mathbb{Z}}
\def\FF{\mathbb{F}}

\begin{document}

\title[Mutation classes of finite type cluster algebras]{Mutation classes of finite type cluster algebras with principal coefficients}

\author{Ahmet I. Seven}

\address{Middle East Technical University, 06531, Ankara, Turkey}
\email{aseven@metu.edu.tr}

\thanks{The author's research was supported in part by Turkish Research Council (TUBITAK)}

\date{June 28th, 2011}

\begin{abstract}
Cluster algebras of finite type is a fundamental class of algebras whose classification is identical to the Cartan-Killing classification. More recently, Fomin and Zelevinsky introduced another central notion of cluster algebras with principal coefficients. These algebras are determined combinatorially by mutation classes of certain rectangular matrices. It was conjectured, by Fomin and Zelevinsky, that finite type cluster algebras with principal coefficients are characterized by the mutation classes which are finite. In this paper, we prove this conjecture.

\end{abstract}

\maketitle

\section{Introduction}
\label{sec:intro}


Cluster algebras of finite type is a fundamental class of algebras whose classification, due to S. Fomin and A. Zelevinsky, is identical to the Cartan-Killing classification \cite{CAII}. More recently, Fomin and Zelevinsky introduced another central notion of cluster algebras with principal coefficients \cite{CAIV}. These algebras are determined combinatorially by mutation classes of certain rectangular matrices. It was conjectured in \cite{CAIV} that finite type cluster algebras with principal coefficients are characterized by the mutation classes which are finite. In this paper we prove this conjecture, using linear algebraic and combinatorial methods.

To be more specific, we need some terminology. In this paper, we deal with the combinatorial aspects of the theory of cluster algebras, so we will not need their definition nor their algebraic properties. The main objects of our study will be skew-symmetrizable matrices and their extensions. Let us recall that an integer matrix $B$ of size $n$ is \emph{skew-symmetrizable} if $DB$ is skew-symmetric for some diagonal matrix $D$ with positive diagonal entries. 
For $m \geq n$, we denote by $\tilde{B}$ a $m\times n$ matrix whose principal, i.e top $n\times n$, part is $B$. Then,
for any matrix index $k$ with $1\leq k \leq n$, the mutation of $\tilde{B}$ at $k$ is the matrix $\mu_k(\tilde{B})=\tilde{B}'$: 
\begin{displaymath}
B' = \left\{ \begin{array}{ll}
B'_{i,j}=-B_{i,j} & \textrm{if $i=k$ or $j=k$}\\
B'_{i,j}=B_{i,j}+sgn(B_{i,k})[B_{i,k}B_{k,j}]_+ & \textrm{else}
\end{array} \right.
\end{displaymath}
(where we use the notation $[x]_+=max\{x,0\}$ and $sgn(x)=x/|x|$ with $sgn(0)=0$). 
Note that the principal part of $\tilde{B'}$ is the matrix $B'=\mu_k(B)$, which is skew-symmetrizable. Mutation is an involutive operation, so repeated mutations give rise to the \emph{mutation-equivalence} relation on $m\times n$ matrices with a skew-symmetrizable principal part. 
The corresponding equivalence classes are called mutation classes. 
A matrix $\tilde{B}$ is said to be of finite mutation type if its mutation class is finite, i.e. only finitely many matrices can be obtained from $\tilde{B}$ by repeated matrix mutations.


Among the matrices $\tilde{B}$, a particular type has turned out to be very special. More precisely, a rectangular matrix associated with principal coefficients is a $2n\times n$ matrix $\tilde{B_0}$ whose top $n\times n$ part $B_0$ is skew-symmetrizable and bottom $n\times n$ part is the identity matrix. Mutation classes of these matrices correspond to cluster algebras with principal coeeficients, which play a central role in the theory of cluster algebras \cite{CAIV}. Another important type of cluster algebras is formed by cluster algebras of finite type. In our set up, these algebras correspond to the mutation classes of skew-symmetrizable matrices of \emph{finite type}; more explicitly, we say that a skew-symmetrizable matrix $B$ is of finite type if, for any $B'$ which is mutation-equivalent to $B$, we have $\left| B'_{i,j}B'_{j,i}\right| \leq 3$. Remarkably, the classification of {finite type} skew-symmetrizable matrices under mutation-equivalence is another instance of Cartan-Killing classification \cite{CAII}. It was conjectured more recently that finite-type skew-symmetrizable matrices are characterized as those whose extended rectangular matrices associated with principal coefficients are of finite mutation type. In this paper we prove this conjecture. More precisely, we obtain the following statement:

\begin{theorem}
\label{th:finite-type-via-2nbyn-matrices}
\cite[Conjecture~4.8]{CAIV}
A skew-symmetrizable $n\times n$ matrix $B_{0}$ is of
finite type if and only if the corresponding $2n \times n$ matrix
$\tilde B_{0}$ is of finite mutation type (where $\tilde B_{0}$ is as defined above).
\end{theorem}

The "only if" part of the conjecture was obtained in \cite{CAIV}; we prove the "if" part, i.e. show that 

($*$) if $B_{0}$ is of infinite type, then $\tilde B_{0}$ is of infinite mutation type. 

\noindent
For this purpose, a convenient setup is provided by a well-known construction that represents skew-symmetrizable matrices by graphs. More precisely, to a skew-symmetrizable matrix $B$ of size $n$, we associate a directed graph $\Gamma(B)$, called the diagram of $B$, with  vertices $1,...,n$ such that there is a directed edge from $i$ to $j$ if and only if $B_{ji}> 0$, and this edge is assigned the weight $|B_{ij}B_{ji}|\,$. Then a mutation $\mu_k$ can be viewed as a transformation on diagrams (see Section~\ref{sec:pre} for a description). To prove ($*$), we restrict ourselves, without losing any generality, to $B_0$ whose diagram is of \emph{minimal infinite type}. These diagrams have been obtained explicitly in \cite{S2} and they are known to be, with few exceptions, mutation-equivalent to the extended Dynkin diagrams (Figure~\ref{fig:extended-dynkin-diagrams}). On the other hand, a description of the mutation classes of extended Dynkin diagrams has been obtained in \cite{S3} using a notion of quasi-Cartan companions, which is a natural generalization of (generalized) Cartan matrices in Kac-Moody Lie algebras \cite{K}. 
This allows us to understand the principal part of $\tilde{B}$ which is mutation-equivalent to $\tilde B_{0}$ in ($*$). However, these methods do not generalize immediately to the whole rectangular matrix $\tilde{B}$. To achieve a generalization, we use an idea of mod $2$ reduction by considering 
$\tilde{B}$ naturally as an alternating bilinear form $\bar{\Omega}^\bullet$ on $\ZZ^{2n}/2\ZZ^{2n}$, which is a vector space over the two-element field. Then, using linear algebraic properties $\bar{\Omega}^\bullet$ along with some properties of generalized Cartan matrices, we reduce ($*$) to the case where $B_0$ has size $2$ and prove the statement.

We prove Theorem~\ref{th:finite-type-via-2nbyn-matrices} in Section~\ref{sec:proof} after some preparation in Section~\ref{sec:pre}.

\section{Preliminaries}
\label{sec:pre}
In this section, we will recall some more terminology and prove some statements that we will use to prove our results. 
First, let us recall that skew-symmetrizable matrices are characterized as follows \cite[Lemma~7.4]{CAII}: $B$ is skew-symmetrizable if and only if $B$ is sign-skew-symmetric (i.e. for any $i,j$ either $B_{i,j}=B_{j,i}=0$ or $B_{i,j}B_{j,i}<0$) and 
for all $k \geq 3$ and all 
$i_1, \dots, i_k\,$, it satisfies
\begin{equation} 
\label{eq:cycle=cycle}
B_{i_1,i_2} B_{i_2,i_3} \cdots B_{i_k,i_1} = 
(-1)^k B_{i_2,i_1} B_{i_3,i_2} \cdots B_{i_1,i_k}\,. 
\end{equation}

%

In this paper, it will be convenient for us to use the following description of the mutation operation on skew-symmetrizable matrices: 
\begin{proposition}
\label{pr:base-change}
Suppose $B$ is a skew-symmetrizable matrix of size $n$ with a skew-symmetrizing matrix $D=diag(d_1,..,d_n)$ (so $d_iB_{i,j}=-d_jB_{j,i}$). Suppose that $\mathcal{B}=\{e_1,...,e_n\}$ is a basis of $\ZZ^n$ and let $\Omega$ be the skew-symmetric bilinear form defined as $\Omega(e_i,e_j)=d_i{B}_{i,j}$, i.e. $DB$ is the Gram matrix of $\Omega$ with respect to the basis $\mathcal{B}$. Then $D\mu_k(B)$ is the Gram matrix of $\Omega$ with respect to the basis $\mathcal{B'}=\{e'_1,e'_2,...,e'_n\}$ defined as follows: $e'_k=-e_k$; $e'_i=e_i-B_{k,i}e_k$ if $B_{k,i}<0$; $e'_i=e_i$ if else. 

\end{proposition}

\proof Note that $\Omega$ is skew-symmetric by its definition $\Omega(e_i,e_j)=d_i{B}_{i,j}=-d_j{B}_{j,i}$. 
Also recall that $B$ and $B'=\mu_k(B)$ shares the same skew-symmetrizing matrix (so $d_iB'_{i,j}=-d_jB'_{j,i}$). Let us now
We will show that $\Omega(e'_i,e'_j)=d_i{B'}_{i,j}$ (which is equal to $=-d_j{B'}_{j,i}$). 
Let us first consider the case with $j=k$. Then, for any $i\ne k$, we have the following: if $B_{k,i}\geq 0$, then $\Omega(e'_i,e'_k)=\Omega(e_i,-e_k)=-d_iB_{i,k}=d_iB'_{i,k}$ 
; similarly, if $B_{k,i}< 0$, then $\Omega(e'_i,e'_k)=\Omega(e_i-B_{k,i}e_k,-e_k)=\Omega(e_i,-e_k)+B_{k,i}\Omega(e_k,e_k)=-\Omega(e_i,e_k)=-d_iB_{i,k}=d_iB'_{i,k}$. Let us now consider the case where $i,j\ne k$. First assume that $B_{k,i}$ and $B_{k,j}$ have the same sign (then $B_{i,j}=B'_{i,j}$), say both are less than zero. Then 
$\Omega(e'_i,e'_j)=\Omega(e_i-B_{k,i}e_k,e_j-B_{k,j}e_k)=\Omega(e_i,e_j)-B_{k,j}\Omega(e_i,e_k)-B_{k,i}\Omega(e_k,e_j)+B_{k,i}^2\Omega(e_k,e_k)=\Omega(e_i,e_j)-B_{k,j}\Omega(e_i,e_k)-B_{k,i}\Omega(e_k,e_j)$. Replacing $\Omega(e_i,e_k)=-d_kB_{k,i}$ and $\Omega(e_k,e_j)=d_k{B}_{k,j}$, we have 
$\Omega(e'_i,e'_j)=\Omega(e_i,e_j)-B_{k,j}(-d_kB_{k,i})-B_{k,i}d_k{B}_{k,j}=\Omega(e_i,e_j)+B_{k,j}(d_kB_{k,i})-B_{k,i}d_k{B}_{k,j}=\Omega(e_i,e_j)$. Now assume that $B_{k,i}$ and $B_{k,j}$ have opposite signs, say $B_{k,i}<0$ and $B_{k,j}>0$ (then $B'_{i,j}=B_{i,j}+B_{i,k}B_{k,j}$ and $B'_{j,i}=B_{j,i}-B_{j,k}B_{k,i}$ ). Then $\Omega(e'_i,e_j')=\Omega(e_i-B_{k,i}e_k,e_j)=\Omega(e_i,e_j)-B_{k,i}\Omega(e_k,e_j)=-d_j{B}_{j,i}-B_{k,i}(-d_j{B}_{j,k})= -d_j({B}_{j,i}-B_{k,i}{B}_{j,k})=-d_j{B'}_{j,i}=d_i{B'}_{i,j}$ as required.

\endproof


In this paper, we also consider mutations of rectangular matrices (Section~\ref{sec:intro}). It is possible to view mutation of a rectangular matrix as a mutation of a skew-symmetrizable matrix as follows:

\begin{definition} 
\label{def:bullet} 

Suppose $m>n$ and $B$ is a skew-symmetrizable $n\times n$ matrix $B$. Let $\tilde{B}$ be a $m\times n$ matrix such that 
$\tilde{B}_{i,j}=B_{i,j}$ (so the top $n\times n$ part of $\tilde{B}$ is $B$). Let $L$ denote the lower $(m-n)\times n$ part of $\tilde{B}$.
We denote by $\tilde{B}^\bullet$ the $m\times m$ matrix which extends $\tilde{B}$ as follows: the left $m\times n$ part is $\tilde{B}$, the upper $n\times (m-n)$ part is $-L$ and the lower-right $(m-n) \times (m-n)$ part is the zero matrix. The matrix $\tilde{B}^\bullet_0$ is skew-symmetrizable.
Furthermore, for $1\leq k \leq n$, the matrix $\mu_k(\tilde{B})$ is the left $m\times n$ part of $\mu_k(\tilde{B}^\bullet)$.

\end{definition} 
\noindent
(Note that the lower-right $(m-n) \times (m-n)$ part of the notation $\tilde{B}^\bullet$ is inconsequential for our study in this paper, so it could have been taken as any matrix.)

\subsection{Diagrams of skew-symmetrizable matrices and their mutations}
\label{subsec:skew-symmetrizable} 
Suppose that $B$ is a skew-symmetrizable matrix of size $n$. Then, following the convention in \cite{DWZ}, the \emph{diagram} of $\tilde{B}$ is the directed graph $\Gamma ({B})$ defined as follows: the vertices of $\Gamma ({B})$ are the indices $1,2,...,n$ such that there is a directed edge from $i$ to $j$ if and only if ${B}_{ji} > 0$, and this edge is assigned the weight $|B_{ij}B_{ji}|\,$. For a rectangular matrix $\tilde{B}$ we define its diagram $\Gamma(\tilde{B})$ as the diagram of the skew-symmetrizable matrix $\tilde{B}^\bullet$ (Definition~\ref{def:bullet}).

It follows from \eqref{eq:cycle=cycle} that the diagram $\Gamma(B)$ of any skew-symmetrizable matrix $B$ has the following property:


\begin{align}
\label{eq:perfect-sq}
&\text{the product of weights along any cycle is a perfect square, i.e. the square}
\\
\nonumber
&\text{of an integer. }
\end{align}

\noindent
Thus we can use the term diagram to mean a directed graph, with no loops or two-cycles, such that the edges are weighted with positive integers satisfying \eqref{eq:perfect-sq}. Let us note that if an edge in a diagram has weight equal to one, then we do not specify its weight in the picture. 


Let us note that if $B$ is not skew-symmetric, then the diagram $\Gamma(B)$ does not determine $B$ as there could be several different skew-symmetrizable matrices whose diagrams are equal (this property will be useful to us in Lemma~\ref{lem:di=ui=1}); however, if a skew-symmetrizing matrix $D$ is fixed, then $\Gamma(B)$ determines $B$.

We also use the following terminology related to diagrams. By a \emph{subdiagram} of $\Gamma$, we always mean a diagram $\Gamma'$ 
obtained from $\Gamma$ by taking an induced (full) directed subgraph on a subset of vertices and keeping all its edge weights the same as in $\Gamma$ 
\cite[Definition~9.1]{CAII}. By a cycle we mean a subdiagram whose vertices can be labeled by elements of $\ZZ/m\ZZ$ so that the edges betweeen them are precisely $\{i,i+1\}$ for $i \in  \ZZ/m\ZZ$. A diagram is called \emph{acyclic} if it has no oriented cycles at all. We call a vertex $v$ \emph{source} (\emph{sink}) if all incident edges are oriented away (towards) $v$. It is well-known that an acyclic diagram has a source and a sink.

\begin{figure}[ht] 
\begin{center}
\setlength{\unitlength}{1.5pt} 
\begin{picture}(30,17)(-5,0) 
\put(0,0){\line(1,0){20}} 
\put(0,0){\line(2,3){10}} 
\put(0,0){\vector(2,3){6}} 
\put(10,15){\line(2,-3){10}} 
\put(10,15){\vector(2,-3){6}} 
\put(0,0){\circle*{2}} 
\put(20,0){\circle*{2}} 
\put(10,15){\circle*{2}} 
\put(2,10){\makebox(0,0){$\alpha$}} 
\put(18,10){\makebox(0,0){$\beta$}} 
\put(10,-4){\makebox(0,0){$\gamma$}} 
\put(10,19){\makebox(0,0){$k$}} 
\end{picture} 
$ 
\begin{array}{c} 
\stackrel{\textstyle\mu_k}{\longleftrightarrow} 
\\[.3in] 
\end{array} 
$ 
\setlength{\unitlength}{1.5pt} 
\begin{picture}(30,17)(-5,0) 
\put(0,0){\line(1,0){20}} 
\put(0,0){\line(2,3){10}} 
\put(10,15){\vector(-2,-3){6}} 
\put(10,15){\line(2,-3){10}} 
\put(20,0){\vector(-2,3){6}} 
\put(0,0){\circle*{2}} 
\put(20,0){\circle*{2}} 
\put(10,15){\circle*{2}} 
\put(2,10){\makebox(0,0){$\alpha$}} 
\put(18,10){\makebox(0,0){$\beta$}} 
\put(10,-4){\makebox(0,0){$\gamma'$}} 
\put(10,19){\makebox(0,0){$k$}} 
\end{picture} 
\end{center}
 
\vspace{-.2in} 
\caption{Diagram mutation} 
\label{fig:diagram-mutation-general} 
\end{figure}
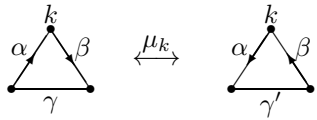

For any vertex $k$ in a diagram $\Gamma$, the associated mutation $\mu_k$ is the transformation that changes $\Gamma$ in such a way that $\mu_k(\Gamma(B))=\Gamma(\mu_k(B))$. More explicitly the mutation $\mu_k$ changes $\Gamma$ as follows \cite{CAII}:
\begin{itemize} 
\item The orientations of all edges incident to~$k$ are reversed, 
their weights intact. 
\item 
For any vertices $i$ and $j$ which are connected in 
$\Gamma$ via a two-edge oriented path going through~$k$ (see  
Figure~\ref{fig:diagram-mutation-general}), 
the direction of the edge $\{i,j\}$ in $\mu_k(\Gamma)$ and its weight $\gamma'$ are uniquely determined by the rule 
\begin{equation} 
\label{eq:weight-relation-general} 
\pm\sqrt {\gamma} \pm\sqrt {\gamma'} = \sqrt {\alpha\beta} \,, 
\end{equation} 
where the sign before $\sqrt {\gamma}$ 
(resp., before $\sqrt {\gamma'}$) 
is ``$+$'' if $i,j,k$ form an oriented cycle 
in~$\Gamma$ (resp., in~$\mu_k(\Gamma)$), and is ``$-$'' otherwise. 
Here either $\gamma$ or $\gamma'$ can be equal to~$0$, which means that the corresponding edge is absent. 
 
\item 
The rest of the edges and their weights in $\Gamma$ 
remain unchanged. 
\end{itemize} 

\noindent 
This operation is involutive, i.e. $\mu_k(\mu_k(\Gamma))=\Gamma$, so it defines an equivalence relation on the set of all diagrams. More precisely, two diagrams are called \emph{mutation-equivalent} if they can be obtained from each other by applying a sequence of mutations. The \emph{mutation class} of a diagram $\Gamma$ is the set of all diagrams which are mutation-equivalent to $\Gamma$.

\subsection{Minimal infinite type diagrams}
\label{subsec:finite}
Recall that a skew-symmetrizable matrix $B$ is said to be of \emph{finite type} if any matrix $B'$ which is mutation-equivalent to $B$ satisfies $|B'_{i,j}B'_{j,i}|\leq 3$ for all $i$ and $j$. We say that $B$ is of \emph{infinite type} if it is not of finite type.
Similarly, a diagram $\Gamma$ is said to be of \emph{finite type} if any diagram $\Gamma'$ 
which is mutation-equivalent to $\Gamma$ has all edge weights equal to $1,2$ or~$3$. 
We say that $\Gamma$ is of \emph{minimal infinite}\footnote{the term minimal 2-infinite was used in \cite{S2}} type if it is of infinite type and any proper subdiagram is of finite type.





Note that $B$ is of finite type if and only if its diagram $\Gamma(B)$ is of finite type.
Let us also note that a subdiagram of a finite type diagram is also of finite type; obviously every diagram which is mutation-equivalent to a diagram of finite type is of finite type itself. Furthermore, a diagram is of finite type if and only if it does not contain any minimal infinite type diagram as a subdiagram.

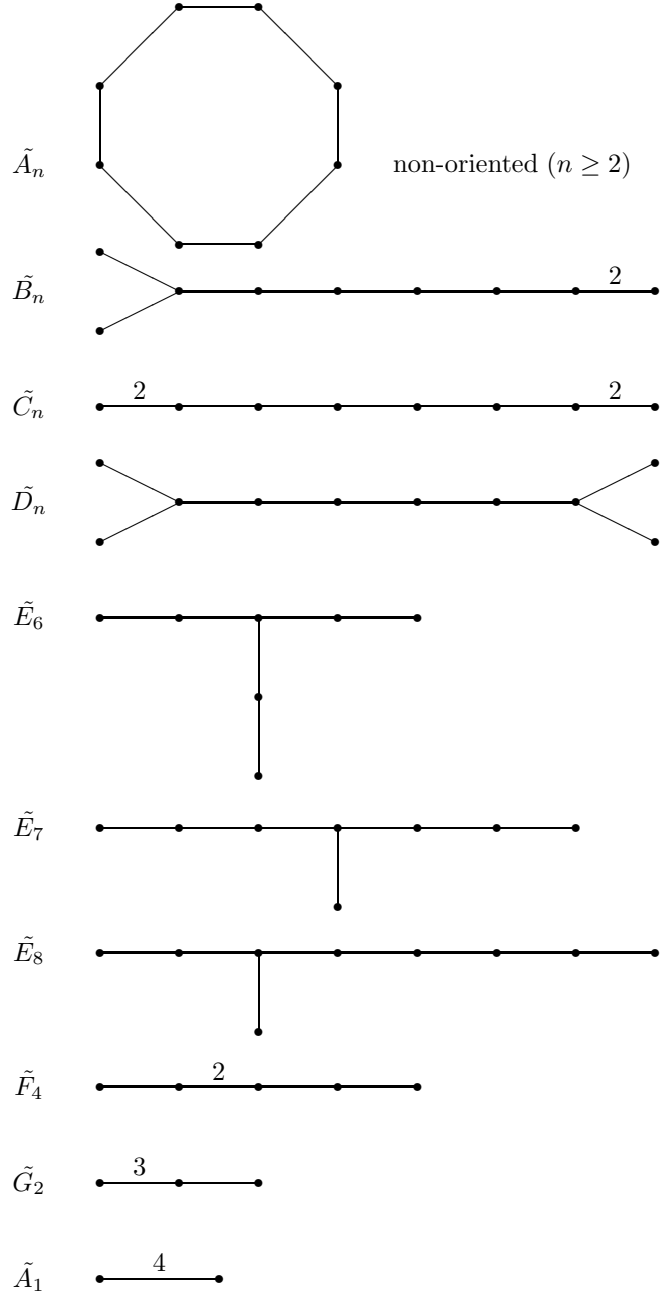
\begin{figure}[ht] 
\[ 
\begin{array}{ccl} 
\tilde{A_n}
&& 
\setlength{\unitlength}{1.5pt} 
\begin{picture}(140,60)(0,-2) 
\put(60,0){\circle*{2.0}} 
\put(60,20){\circle*{2.0}}
\put(40,40){\circle*{2.0}}
\put(20,40){\circle*{2.0}}
\put(0,20){\circle*{2.0}}
\put(20,-20){\circle*{2.0}}
\put(0,0){\circle*{2.0}}
\put(40,-20){\circle*{2.0}}

\put(60,0){\line(-1,-1){20}}
\put(60,0){\line(0,1){20}}
\put(60,20){\line(-1,1){20}}
\put(40,40){\line(-1,0){20}}
\put(20,40){\line(-1,-1){20}}
\put(0,20){\line(0,-1){20}}
\put(0,0){\line(1,-1){20}}
\put(20,-20){\line(1,0){20}}

\put(104,0){\makebox(0,0){non-oriented ($n\geq 2$)}}
\end{picture} 
\\[.3in] 
\tilde{B_n}
&& 
\setlength{\unitlength}{1.5pt} 
\begin{picture}(140,15)(0,-2) 
\put(20,0){\line(1,0){120}} 
\put(0,10){\line(2,-1){20}} 
\put(0,-10){\line(2,1){20}} 
\multiput(20,0)(20,0){7}{\circle*{2}} 
\put(0,10){\circle*{2}} 
\put(0,-10){\circle*{2}} 
\put(130,4){\makebox(0,0){$2$}} 
\end{picture} 
\\[.2in] 
\tilde{C_n}
&& 
\setlength{\unitlength}{1.5pt} 
\begin{picture}(140,17)(0,-2) 
\put(0,0){\line(1,0){140}} 
\multiput(0,0)(20,0){8}{\circle*{2}} 
\put(10,4){\makebox(0,0){$2$}} 
\put(130,4){\makebox(0,0){$2$}} 
\end{picture} 
\\[.1in] 
\tilde{D_n}
&& 
\setlength{\unitlength}{1.5pt} 
\begin{picture}(140,17)(0,-2) 
\put(20,0){\line(1,0){100}} 
\put(0,10){\line(2,-1){20}} 
\put(0,-10){\line(2,1){20}} 
\put(120,0){\line(2,-1){20}} 
\put(120,0){\line(2,1){20}} 
\multiput(20,0)(20,0){6}{\circle*{2}} 
\put(0,10){\circle*{2}} 
\put(0,-10){\circle*{2}} 
\put(140,10){\circle*{2}} 
\put(140,-10){\circle*{2}} 
\end{picture} 
\\[.2in] 
\tilde{E_6}
&& 
\setlength{\unitlength}{1.5pt} 
\begin{picture}(140,17)(0,-2) 
\put(0,0){\line(1,0){80}} 
\put(40,0){\line(0,-1){40}} 
\put(40,-20){\circle*{2}} 
\put(40,-40){\circle*{2}} 
\multiput(0,0)(20,0){5}{\circle*{2}} 
\end{picture} 
\\[.7in] 
\tilde{E_7}
&& 
\setlength{\unitlength}{1.5pt} 
\begin{picture}(140,17)(0,-2) 
\put(0,0){\line(1,0){120}} 
\put(60,0){\line(0,-1){20}} 
\put(60,-20){\circle*{2}} 
\multiput(0,0)(20,0){7}{\circle*{2}} 
\end{picture} 
\\[.25in] 
\tilde{E_8}
&& 
\setlength{\unitlength}{1.5pt} 
\begin{picture}(140,17)(0,-2) 
\put(0,0){\line(1,0){140}} 
\put(40,0){\line(0,-1){20}} 
\put(40,-20){\circle*{2}} 
\multiput(0,0)(20,0){8}{\circle*{2}} 
\end{picture} 
\\[.3in] 
\tilde{F_4}
&& 
\setlength{\unitlength}{1.5pt} 
\begin{picture}(140,17)(0,-2) 
\put(0,0){\line(1,0){80}} 
\multiput(0,0)(20,0){5}{\circle*{2}} 
\put(30,4){\makebox(0,0){$2$}} 
\end{picture} 
\\[.1in] 
\tilde{G_2} 
&& 
\setlength{\unitlength}{1.5pt} 
\begin{picture}(140,17)(0,-2) 
\put(0,0){\line(1,0){40}} 
\multiput(0,0)(20,0){3}{\circle*{2}} 
\put(10,4){\makebox(0,0){$3$}} 
\end{picture}\\[.1in] 
\tilde{A_1}
&& 
\setlength{\unitlength}{1.5pt} 
\begin{picture}(65,17)(20,-2) 
\put(20,0){\line(1,0){30}} 
\put(20,0){\circle*{2}} 
\put(50,0){\circle*{2}} 
\put(35,4){\makebox(0,0){$4$}} 
\end{picture} 
\end{array} 
\] 
\caption{Extended Dynkin diagrams are orientations of the extended Dynkin graphs given above; the first graph $\tilde{A_n}$ 
is assumed to be a  non-oriented cycle, the rest of the graphs are assumed to be arbitrarily oriented; each $\tilde{X_n}$ 
has $n+1$ vertices} 
\label{fig:extended-dynkin-diagrams} 
\end{figure}

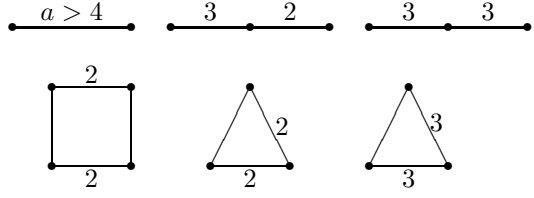
\begin{figure}[ht]

\begin{center}
\setlength{\unitlength}{1.5pt}
\begin{picture}(130,44)(0,12)

\put(10,0){\circle*{2.0}}
\put(10,20){\circle*{2.0}}
\put(30,20){\circle*{2.0}}
\put(30,0){\circle*{2.0}}

\put(50,0){\circle*{2.0}}
\put(70,0){\circle*{2.0}}
\put(60,20){\circle*{2.0}}
\put(90,0){\circle*{2.0}}
\put(110,0){\circle*{2.0}}
\put(100,20){\circle*{2.0}}

\put(107,11){\makebox(0,0){$3$}}

\put(20,-3){\makebox(0,0){$2$}}
\put(20,23){\makebox(0,0){$2$}}

\put(60,-3){\makebox(0,0){$2$}}
\put(68,10){\makebox(0,0){$2$}}
\put(100,-3){\makebox(0,0){$3$}}

\put(10,0){\line(1,0){20}}
\put(10,0){\line(0,1){20}}
\put(10,20){\line(1,0){20}}
\put(30,0){\line(0,1){20}}

\put(50,0){\line(1,0){20}}
\put(60,20){\line(1,-2){10}}
\put(60,20){\line(-1,-2){10}}

\put(90,0){\line(1,0){20}}
\put(90,0){\line(1,2){10}}
\put(100,20){\line(1,-2){10}}

\put(0,35){\line(1,0){30}} 
\put(0,35){\circle*{2}} 
\put(30,35){\circle*{2}} 
\put(15,39){\makebox(0,0){$a>4$}} 
\put(40,35){\line(1,0){40}} 
\put(50,39){\makebox(0,0){$3$}} 
\put(40,35){\circle*{2}} 
\put(60,35){\circle*{2}} 
\put(80,35){\circle*{2}} 
\put(70,39){\makebox(0,0){$2$}} 

\put(90,35){\line(1,0){40}} 
\put(100,39){\makebox(0,0){$3$}} 
\put(90,35){\circle*{2}} 
\put(110,35){\circle*{2}} 
\put(130,35){\circle*{2}} 
\put(120,39){\makebox(0,0){$3$}} 


\end{picture}\\[.6in]
\end{center}

\caption{Minimal infinite type diagrams which are not mutation equivalent to any extended Dynkin diagram: each diagram has an arbitrary acyclic orientation} 
\label{fig:minimal-non-extended-Dynk} 
\end{figure}

Classification of finite type diagrams, by Fomin and Zelevinsky, is another instance of Cartan-Killing classification \cite{CAII}:
a skew-symmetrizable matrix $B$ with a connected diagram $\Gamma(B)$ is 
of finite type if and only if $\Gamma(B)$ is mutation-equivalent to an arbitrarily oriented Dynkin diagram. 
An alternative characterization was obtained in \cite{S2} by giving a complete list of the minimal infinite type diagrams. In this paper, we will use the following property of these diagrams: any minimal infinite type diagram is either mutation-equivalent to an extended Dynkin diagram (Figure~\ref{fig:extended-dynkin-diagrams}) or it is one of the diagrams in Figure~\ref{fig:minimal-non-extended-Dynk}. 



\subsection{Symmetrizable matrices and their mutations}
\label{subsec:symmetrizable}
Let us now recall "quasi-Cartan companions" that we use to study mutation classes of extended Dynkin diagrams.


\begin{definition} 
\label{def:symmetrizable} 
Let $A$ be a $n \times n$ matrix (whose entries are integers). 
The matrix $A$ is called symmetrizable if there exists a diagonal matrix $D$ with positive diagonal entries such that $DA$ is symmetric. We say that $A$ is a quasi-Cartan matrix if it is symmetrizable and all of its diagonal entries are equal to $2$.
\end{definition}
\noindent
The symmetrizable matrix $A$ is sign-symmetric, i.e. $sgn(A_{i,j})=sgn(A_{j,i})$. We say that $A$ is (semi)positive if $DA$ is positive (semi)definite, i.e. (resp. $x^TDAx\geq 0$) $x^TDAx>0$ for all $x\ne 0$ (here $x^T$ is the transpose of $x$ which is a vector viewed as a column matrix). We say that $u$ is a \emph{radical} vector of $A$ if $Au=0$; we call $u$ \emph{sincere} if all of its coordinates are non-zero. 
A quasi-Cartan matrix is a \emph{generalized Cartan matrix} if all of its non-zero entries which are not on the diagonal are negative.  


Quasi-Cartan matrices are related to skew-symmetrizable matrices via the following notion:
\begin{definition} 
\label{def:companion} 
Let $B$ be a skew-symmetrizable matrix. A \emph{quasi-Cartan companion} (or "companion" for short) of $B$ is a quasi-Cartan matrix $A$ with $|A_{i,j}|= |B_{i,j}|$ for all $i \ne j$. We say that $A$ is \emph{admissible} if it satisfies the following sign condition: for any cycle $Z$ in $\Gamma$, the product $\prod _{\{i,j\}\in Z}(-A_{i,j})$ over all edges of $Z$ is negative if $Z$ is oriented and positive if $Z$ is non-oriented.


\end{definition}
\noindent
Let us note that a skew-symmetrizing matrix $D$ for $B$ is a symmetrizing matrix for a quasi-Cartan companion $A$.  
Let us also note that we may view a quasi-Cartan companion $A$ of $B$ as a sign assignment to the edges (of the underlying undirected graph) of $\Gamma=\Gamma(B)$; more explicitly an edge $\{i,j\}$ is assigned the sign of the entry $A_{i,j}$ (which is the same as the sign of $A_{j,i}$ because $A$ is sign-symmetric). Then
the sign condition in the definition can also be described as follows: if $Z$ is (non)oriented, then there is exactly an (resp. even) odd number of edges $\{i,j\}$ such that $(A_{i,j})>0$. (recall that, since $A$ is symmetrizable, we have $sgn(A_{i,j})=sgn(A_{j,i})$). 
Any two admissible companions of $B$ can be obtained from each other by a sequence of simultaneous sign changes in rows and columns (see \cite[Section~2]{S3} for more details).


To be able to use symmetrizable matrices to study the mutation operation, we use the following extension of the mutation operation to quasi-Cartan companions \cite[Proposition~3.2]{BGZ}:

\begin{definition} 
\label{def:comp-mut} 
Suppose that $\Gamma$ is a diagram and let $A$ be a quasi-Cartan companion of $\Gamma$. 
Let $k$ be a vertex in $\Gamma$. "The mutation of $A$ at $k$" is the quasi-Cartan matrix $A'$ such that for any $i,j \ne k$: $A'_{i,k}=sgn(B_{i,k})A_{i,k}$, $A'_{k,j}=-sgn(B_{k,j})A_{k,j}$, $A'_{i,j}=A_{i,j}-sgn(A_{i,k}A_{k,j})[B_{i,k}B_{k,j}]_+$. 
It is a quasi-Cartan companion of $\mu_k(\Gamma)$ if $A$ is admissible \cite[Proposition~3.2]{BGZ}.
\end{definition}

\noindent
This operation may also be viewed as a base change for a symmetric bilinear form in a way similar to the one for skew-symetrizable matrices as in Proposition~\ref{pr:base-change}. To be more precise, suppose that $D$ is a skew-symetrizing matrix of $B$. Then $D$ is also a symmetrizer for $A$, with $DA=C$ symmetric. If we consider $C$ as the Gram matrix of a symmetric bilinear form on $\ZZ^n$ with respect to a basis $\mathcal{A}=\{e_1,...,e_n\}$, then $DA'=C'$ is the Gram matrix of the same symmetric bilinear form with respect to the basis $\mathcal{A'}=\{e'_1,e'_2,...,e'_n\}$ defined as follows: $e'_k=-e_k$; $e'_i=e_i-A_{k,i}e_k$ if $B_{k,i}<0$; $e'_i=e_i$ if else. Note that if $(u'_1,...,u'_n)$ are the coordinates of $u=(u_1,..,u_n)$ with respect to the new basis $\mathcal{A'}$, then we have the following: for $i\ne k$, $u'_i=u_i$;   $u'_k=-u_k+\sum A_{k,i}u_i$ where the sum is over $i$ with $B_{k,i}<0$. Let us also note that the basis $\mathcal{A'}$ coincides with $\mathcal{B'}$ of Proposition~\ref{pr:base-change}
 in $\ZZ^n/2\ZZ^n$ as a vector space over the field $\FF_2$ with two elements.

Let us also note that $A'$ may not be admissible even if $A$ is admissible. However admissibility is preserved for the matrices we are interested in this paper:


\begin{proposition}
\label{pr:admissible-preserved-extended-Dynkin}\cite[Proposition~5.1]{S3}
Suppose that $B$ is a skew-symmetrizable matriz whose diagram $\Gamma(B)$ is mutation-equivalent to an extended Dynkin diagram. Let $A$ be an admissible quasi-Cartan companion of $B$ and let $A'$ be the mutation of $A$ at $k$. Then $A'$ is an admissible quasi-Cartan companion of $\mu_k(B)$. Furthermore $A$ is semipositive with corank $1$.

\end{proposition}

Among the diagrams which are mutation-equivalent to an extended Dynkin diagram, minimal infinite diagrams are characterized as those with an admissible companion whose radical vector is sincere \cite[Theorem~3.4]{S3}. Also, a diagram is of finite type if and only if it has an admissible quasi-Cartan companion which is positive \cite[Theorem~1.2]{BGZ}.



Let us now determine some basic types of rectangular matrices which are of infinite mutation type:

\begin{proposition}
\label{pr:m>4-infinite}
Suppose that $B$ is a $2\times 2$ skew-symmetrizable matrix such that $\Gamma(B)$ is a connected two vertex diagram whose edge-weight is 
greater than $4$. Then any rectangular $m\times 2$ matrix $\tilde{B}$, $m\geq 3$, with a connected diagram $\Gamma(\tilde{B})$, is of infinite mutation type.

\end{proposition}

\proof It is enough to prove this for $m=3$. We denote by $\Gamma=\Gamma(B)$ the diagram of $B$, so $\Gamma$ is a two-vertex diagram whose vertices are labeled by  $1$ and $2$. We denote the weight of $\Gamma$ by $\gamma >4$. The diagram $\tilde{\Gamma}=\Gamma(\tilde{B})$ contains $\Gamma$ and has one extra vertex labeled by $3$. 


First suppose that the vertex $3$ is connected to both vertices in $\Gamma$ such that $\tilde{\Gamma}$ is an oriented triangle. Let us first assume that the weight of the edges $\{1,3\}$ and $\{2,3\}$ are $\alpha$ and $\beta$ respectively. Without loss of generality we also assume that $\alpha \leq \beta$. Then $\tilde{\Gamma}'=\mu_2(\tilde{\Gamma})$ is oriented triangle with weights with weights $\beta$, $\gamma$ and $\alpha'=\alpha+\gamma \beta-2\sqrt{\alpha \beta \gamma}$ (this is equal to $(\sqrt{\gamma \beta} - \sqrt{\alpha})^2$). We claim that $\alpha'>\alpha$: Suppose to the contrary that $\alpha'\leq \alpha$. Then $\gamma \beta \leq 2\sqrt{\alpha\beta\gamma}$. Taking squares of both sides gives $\gamma^2\beta^2\leq 4\alpha\beta\gamma$ implying $\gamma\beta\leq 4\alpha$, which contradicts our assumption that $\gamma >4$. Thus the sum of the weights of $\tilde{\Gamma}'$ is greater than the sum of the weights of $\tilde{\Gamma}$. Continuing with applying mutations at the vertices incident to the edge with larger weight, thus increasing the sum of the weights, we see that the mutation class of $\tilde{\Gamma}$ is infinite.

Now assume that the vertex $3$ is connected to $\Gamma$ such that $\tilde{\Gamma}$ is acyclic. Mutating at a vertex in $\Gamma$ if necessary, we may assume that there is a vertex $k$ in $\Gamma$ which is not a source nor a sink in  $\tilde{\Gamma}$. Then  $\mu_k(\tilde{\Gamma})$ is an oriented triangle (containing $\mu_k({\Gamma})$ as a subdiagram having the same weight $\gamma$), so we may apply our previous argument to see that the mutation class of $\tilde{\Gamma}$ is infinite. \endproof

\begin{proposition}
\label{pr:4-infinite}
Suppose that $B$ is a $2\times 2$ skew-symmetrizable matrix such that $\Gamma=\Gamma(B)$ is a connected two vertex diagram whose edge-weight is equal to $4$. Then a rectangular $m\times 2$ matrix $\tilde{B}$, $m\geq 3$, with a connected diagram $\Gamma(\tilde{B})$ is of finite mutation type if and only if, for all $i=3,...,m$, the diagram $T=\{1,2,i\}$ is an oriented triangle such that the weights of the edges $\{i,1\}$ and $\{i,2\}$ are equal.

\end{proposition}

\proof It is enough to prove this for $m=3$. First assume that $T$ is an oriented triangle with the edges $\{i,1\}$ and $\{i,2\}$ having equal weight $\alpha$. Then it follows from a direct check that both $\mu_1(\tilde{\Gamma})$ and $\mu_2(\tilde{\Gamma})$ are oriented triangles with the same weights as in $\tilde{\Gamma}$, thus the mutation class of $\tilde{\Gamma}$ is finite.

Let us now assume that $T$ is oriented with the edges $\{i,1\}$ and $\{i,2\}$ having non-equal weights say $\alpha$ and $\beta$ respectively with $\beta > \alpha$. Then $\widetilde{\Gamma'}=\mu_2(\tilde{\Gamma})$ is oriented triangle with weights $4$, $\beta$ and $\alpha'=\alpha+4 \beta-2\sqrt{4\alpha \beta }$ (this is equal to $(\sqrt{4 \beta} - \sqrt{\alpha})^2$). We claim that $\alpha'>\alpha$; suppose to the contrary that $\alpha'\leq \alpha$. Then $4 \beta \leq 2\sqrt{4\alpha\beta}$. Taking squares of both sides gives $4^2\beta^2\leq 4^2\alpha\beta$ implying $\beta\leq \alpha$, which contradicts our assumption. Thus the sum of the weights of $\widetilde{\Gamma'}$ is greater than the sum of the weights of $\tilde{\Gamma}$. Continuing with applying mutations at the vertices of $\Gamma$ incident to the edge with larger weight, we see that the mutation class of $\tilde{\Gamma}$ is infinite.

Now assume that $T$ is acyclic. Mutating at a vertex in $\Gamma$ if necessary, we may assume that there is a vertex $k$ in $\Gamma$ which is not a source nor a sink in  $T$. Then $\mu_k(T)$ is an oriented triangle with the edges $\{i,1\}$ and $\{i,2\}$ having non-equal weights (with $\mu_k({\Gamma})$ having the same weight $4$), so we may apply our previous argument to see that the mutation class of $\tilde{\Gamma}$ is infinite. This completes the proof of the proposition. \endproof

\newpage

\clearpage

\section{Proof of the main result}
\label{sec:proof}
First we will prove two lemmas that we use to prove Theorem~\ref{th:finite-type-via-2nbyn-matrices}. 


\begin{lemma}
\label{lem:min-inf-non-ext>4}
Suppose that $\Gamma$ is a minimal infinite type diagram which is not mutation-equivalent to any extended Dynkin diagram, i.e. $\Gamma$ is one of the diagrams in Figure~\ref{fig:minimal-non-extended-Dynk}. Then $\Gamma$ is mutation-equivalent to a diagram which contains an edge whose weight is greater than 4.
\end{lemma}

\proof This is obvious if $\Gamma$ has exactly two vertices. If $\Gamma(B_0)$ has exactly three vertices, then mutating at a vertex if necessary, we may assume that there is a vertex $k$ in $\Gamma$ which is not a source nor a sink. Then $\mu_k(\Gamma)$ is an oriented triangle which contains an edge whose weight is greater than $4$. Similarly, if $\Gamma(B_0)$ has four vertices, then mutating at a vertex if necessary, we may assume that there is a vertex $k$ in $\Gamma$ which is not a source nor a sink. Then $\mu_k(\Gamma)$ contains a three-vertex subdiagram which belongs to Figure~\ref{fig:minimal-non-extended-Dynk}, so the statement follows from our previous argument.
\endproof

\begin{lemma}
\label{lem:di=ui=1}
Suppose that $\Gamma$ is a diagram mutation equivalent to an extended diagram which is not of type $\tilde{C}$. Then there is a skew-symmetrizable matrix $B$ with $\Gamma=\Gamma(B)$ such that the following holds: $B$ has a skew-symetrizing matrix $D=diag(d_1,...,d_n)$ and an admissible quasi-Cartan companion $A$ with a non-zero radical vector $u=(u_1,...,u_n)$ such that for some index $l$ we have $d_l=1$ mod 2 and $u_l=1$ mod 2. 

Furthermore if this statement is true for $\Gamma$ with $B,A$ and $u$, then, for any vertex $k$, it is true for $\mu_k(\Gamma)$ with $\mu_k(B), \mu_k(A)$ and $u'$ as in Definition~\ref{def:comp-mut}.

\end{lemma}

\proof If $\Gamma$ is an extended Dynkin diagram which is not of type $\tilde{C}$, then the statement follows from a direct check on tables in \cite[Chapter 4]{K}. More explicitly; for $\tilde{A},\tilde{E_6},\tilde{E_7},\tilde{E_8}$, $\tilde{A_1}$ and $\tilde{G_2}$, we choose a generalized Cartan matrix from Table Aff 1; for $\tilde{B_n}$ and $\tilde{F_4}$ choose from Table Aff 2 in \cite[Chapter 4]{K}. 


Let us now assume that the statement is true for a diagram $\Gamma=\Gamma(B)$ which is mutation equivalent to an extended Dynkin diagram with $B$, $A$ and $u$. We will show that the statement holds for $\Gamma'=\mu_k(\Gamma)=\Gamma(\mu_k(B))$ with $\mu_k(B), \mu_k(A)$ and $u'$. 
We will do this viewing $\mu_k$ as a base change for the associated symmetric bilinear form of $A$ as in Definition~\ref{def:comp-mut}. 
For this let us first recall that $D$ is also a skew-symmetrizing matrix for $\mu_k(B)$ and $\mu_k(A)$. Let $u_1',...,u_n'$ be the coordinates of the $u$ with respect to the new basis (Definition~\ref{def:comp-mut}). Let us note that if $k\ne l$, then $u_l'=u_l$, so the conclusion of the lemma is satisfied. Assume now that $k=l$ and suppose that for any $i\ne k$, we have $d_i=0$ mod 2 or $u_i=0$ mod $2$ (otherwise by the previous argument the lemma is obtained). Let us note that $d_iA_{i,l}=d_lA_{l,i}$ (recall that $B_{i,j}=A_{i,j}$ mod $2$) with $d_l=1$ mod 2, therefore if  $d_i=0$ mod 2 then $A_{l,i}=0$ mod $2$. Thus, under our assumptions, we have either $u_i=0$ mod $2$ or $A_{l,i}=0$ mod 2. Then $u'_l=-u_l+\sum A_{l,i}u_i=u_l=1$ mod $2$ (where the sum is over $i$ with $B_{l,i}<0$). This completes the proof of the lemma. 
\endproof

Let us now prove Theorem~\ref{th:finite-type-via-2nbyn-matrices}. The "only if" part follows from \cite[Proposition~4.9]{CAIV}. We will prove the if part, i.e. we will show that 


%
($*$) if $B_0$ is infinite type then the corresponding $\tilde{B_0}$ is of infinite mutation type. 

\noindent
For this purpose it is enough to show ($*$) 
for $B_0$ such that $\Gamma(B_0)$ is of minimal infinite type. Thus we assume that $\Gamma(B_0)$ is of minimal infinite type with $n$ vertices. If $\Gamma(B_0)$ is not mutation-equivalent to any extended Dynkin diagram, then, by Lemma~\ref{lem:min-inf-non-ext>4}, the diagram $\Gamma(B_0)$ is mutation-equivalent to a diagram which contains an edge whose weight is greater than $4$. Let $\tilde{B}$ be the matrix obtained from $\tilde{B_0}$ by the same sequence of mutations. Then, by Proposition~\ref{pr:m>4-infinite}, the matrix $\tilde{B}$, so $\tilde{B_0}$, is of infinite mutation type.

Let us now suppose that $\Gamma(B_0)$ is mutation-equivalent to an extended Dynkin diagram. 
We denote by $D_0=diag(d_1,...,d_n)$ a skew-symmetrizing matrix of $B_0$ with $S_0=D_0B_0$ is skew-symmetric. 
We also denote by $D=diag(d_1,...,d_n,...,d_{2n})$ a skew-symmetrizing matrix which extends $D_0$ so that $S=D\tilde{B}^\bullet_0$ is skew-symmetric 
(so $d_i \tilde{B}^\bullet_{0_{i,j}}=-d_j\tilde{B}^\bullet_{0_{j,i}}$).

Let $e_i$, $i=1,...,2n$, be the standard basis vectors for $\ZZ^{2n}$. Let $\Omega^\bullet$ be the skew-symmetric bilinear form on $\ZZ^{2n}$ defined by $S=D\tilde{B}^\bullet$ as in Proposition~\ref{pr:base-change} (so $\Omega^\bullet(e_i,e_j)=d_i\tilde{B}^\bullet_{0_{i,j}}$ for $i,j=1,...,2n$). 
We denote by $\bar{\Omega}^\bullet$ the induced alternating bilinear form on $\ZZ^{2n}/2\ZZ^{2n}$ considered as a vector space over the two-element field $\FF_2$. We denote by $\Omega$ the restriction of $\Omega^\bullet$ to $span(e_1,...,e_n)$ (so $\Omega$ is the skew-symmetric form corresponding to $S_0=D_0B_0$). 

We denote by $A_0$ an admissible quasi-Cartan companion of $B_0$ and by $u=(u_1,...,u_n,0,...,0)$ a non-zero radical vector for $\Omega$. 
Note that $u$ is a radical vector for $\bar{\Omega}$ (recall that $B_{0_{i,j}}=A_{0_{i,j}}$ mod $2$). Let us also note that the matrix $A_0$ can be obtained from an affine type generalized Cartan matrix by a sequence of mutations and simultaneous sign changes in rows and columns (Proposition~\ref{pr:admissible-preserved-extended-Dynkin}).


To proceed, we first suppose that $\Gamma(B_0)$ is not mutation equivalent to an extended Dynkin diagram of type $\tilde{C}$ (Figure~\ref{fig:extended-dynkin-diagrams}).
Then, by Lemma~\ref{lem:di=ui=1}, we may assume that for some index $l$, $1\leq l\leq n$, we have  $d_l=1$ mod 2 and $u_l=1$ mod 2. 
Let $r=n+l$. Since $\tilde{B}^\bullet_{0_{r,l}}=1$ and $\tilde{B}^\bullet_{0_{l,r}}=-1$ (by definition of $\tilde{B_0}$ and $\tilde{B}_0^\bullet$)
we have $d_r=d_r\tilde{B}^\bullet_{0_{r,l}}=-d_l\tilde{B}^\bullet_{0_{l,r}}=d_l$,  
in particular $d_r=1$ mod $2$ (recall $d_l=1$ mod 2 by our assumption). Thus $\bar{\Omega}^\bullet(e_{r},e_l)=d_r\tilde{B}^\bullet_{0_{r,l}}=1$ mod $2$. On the other hand, for any $i\ne l$, we have  $\tilde{B}^\bullet_{0_{r,i}}=0$  so $\bar{\Omega}^\bullet(e_{r},e_i)=0$, implying 


\begin{align}
\label{eq:(*)}
&\text{ $\bar{\Omega}^\bullet(e_{r},u)=1$.} 
\end{align}

\noindent
Thus, the vector $u$ is not a radical vector for $\bar{\Omega}^\bullet$, however it is a \emph{non-zero} radical vector for $\bar{\Omega}$, i.e. $\bar{\Omega}(v,u)=0$ for all $v$ in $span(e_1,...,e_n)$. 

Since $\Gamma(B_0)$ is mutation-equivalent to an extended Dynkin diagram, the matrix $B_0$ is mutation-equivalent to a matrix $B=\mu_k...\mu_1(B_0)$ such that $\Gamma(B)$ contains an edge say $\{i,j\}$, ($1\leq i,j\leq n$), whose weight is equal to $4$ \cite{S3}. Then $\tilde{B}^\bullet=\mu_k...\mu_1(\tilde{B}^\bullet_0)$ is the Gram matrix of $\Omega^\bullet$ with respect to the basis $x_1,...,x_{2n}$ obtained from $e_i$ as described in Proposition~\ref{pr:base-change}. Then, in particular,  $span(e_1,...,e_n)=span(x_1,...,x_n)$ and $B$ is the Gram matrix of $\Omega$ with respect to $\mathcal{B}=\{x_1,...,x_n\}$. Furthermore 

\begin{align}
\label{eq:(**)}
&\text{for all $i=1,...,2n$, $x_i=e_i+v_i$ for some $v_i$ in $span(e_1,...,e_n)$.} 
\end{align}

Let us now consider $A=\mu_k...\mu_1(A_0)$ which is the corresponding admissible quasi-Cartan companion of $B$ obtained by the same sequence of mutations (Definition~\ref{def:comp-mut}). Similarly $A$ is the Gram matrix of a symmetric bilinear form with respect to a basis $\mathcal{A}=\{x'_1,...,x'_{n}\}$ of $span(e_1,...,e_n)$ obtained as described in Definition~\ref{def:comp-mut}. Note that, in general $x_i\ne x'_i$, however $x'_i=x_i$ in $\ZZ^n/2\ZZ^n$ for $i=1,...,n$. 


Let us now denote by $A|_{i,j}$ the matrix obtained from $A$ by removing all rows and columns corresponding to the indices which are not in $\{i,j\}$.
Note then that $diag(d_i,d_j)$ is a symmetrizer of $A|_{i,j}$. Also $A|_{i,j}$ is not positive and has corank equal to $1$, therefore $u$ is in the $span(x'_i,x'_j)$ (\cite[Proposition~4.2]{S3}).
Applying a sign change simultaneously at row and a column if necessary, we may assume that $A|_{i,j}$ is a generalized Cartan matrix of size $2$. There are only two types of such matrices: they are $A_1^{(1)}$ and $A_2^{(2)}$ from tables in \cite[Chapter 4]{K}. We note that $A_2^{(2)}$ does not satisfy the conclusion of Lemma~\ref{lem:di=ui=1}; 
however, $A$ satisfies the lemma (by the second part of Lemma~\ref{lem:di=ui=1}), therefore $A|_{i,j}$ is of type $A_1^{(1)}$. This implies that $d_i=d_j$ (also $d_i=d_j=1$ mod 2 because Lemma~\ref{lem:di=ui=1} is satisfied), so $|A_{i,j}|=|A_{j,i}|=2$ (so $|B_{i,j}|=|B_{j,i}|=2$ as well).  
In fact, since we assumed that $A|_{i,j}$ is a generalized Cartan matrix, we have $A_{i,j}=-2$, so then $u=x'_i+x'_j$. Note that $u=x_i+x_j$ in $\ZZ^n/2\ZZ^n$.



Now let us note that by \eqref{eq:(**)} and \eqref{eq:(*)} we have 
$\bar{\Omega}^\bullet(x_r,u)=\bar{\Omega}^\bullet(e_r+v_r,u)=\bar{\Omega}^\bullet(e_r,u)=1$ (because $\bar{\Omega}(v,u)=0$ for all $v$ in $span(e_1,...,e_n)$). Thus ${\Omega}^\bullet(x_r,u)=1$ mod 2. 
On the other hand, ${\Omega}^\bullet(x_r,u)={\Omega}^\bullet(x_r,x_i+x_j)={\Omega}^\bullet(x_r,x_i)+{\Omega}^\bullet(x_r,x_j)=
d_r{\tilde{B}}^\bullet_{r,i}+d_r{\tilde{B}}^\bullet_{r,j}$. Then, since $d_r=1$ mod 2, we have 


\begin{align}
\label{eq:(+)}
&\text{$\bar{\Omega}^\bullet(x_r,u)=\tilde{B}^\bullet_{r,i} + \tilde{B}^\bullet_{r,j}=1$ mod $2$.}
\end{align}

Suppose now, to the contrary of ($*$), 
that $\tilde{B}$ is of finite mutation type. Then, in the diagram $\Gamma(\tilde{B})$ , the triangle $\{r,i,j\}$ is oriented with edges $\{r,i\}$ and $\{r,j\}$ having equal weight (Proposition~\ref{pr:4-infinite}). 
Since $\tilde{B}^\bullet$ is skew-symetrizable, by \eqref{eq:cycle=cycle}, we have $\tilde{B}^\bullet_{j,i}\tilde{B}^\bullet_{r,j}\tilde{B}^\bullet_{i,r}=-\tilde{B}^\bullet_{i,j}
\tilde{B}^\bullet_{j,r}\tilde{B}^\bullet_{r,i}$. Then, since $|B_{i,j}|=|B_{j,i}|=2$, we have

\begin{align}
\label{eq:(***)}
&\text{$2\tilde{B}^\bullet_{r,j}\tilde{B}^\bullet_{i,r}=2
\tilde{B}^\bullet_{j,r}\tilde{B}^\bullet_{r,i}$ }
\end{align}
 
\noindent
Also, since the weights of the edges $\{r,i\}$ and $\{r,j\}$ are equal, we have 

\begin{align}
\label{eq:(****)}
&\text{$\tilde{B}^\bullet_{r,i}\tilde{B}^\bullet_{i,r}=\tilde{B}^\bullet_{j,r}\tilde{B}^\bullet_{r,j}$.} 
\end{align}

\noindent
Now multiplying both sides of \eqref{eq:(****)} by $\tilde{B}^\bullet_{r,j}$, we have

\begin{align}
\label{eq:(v*)}
&\text{$\tilde{B}^\bullet_{r,j}(\tilde{B}^\bullet_{r,i}\tilde{B}^\bullet_{i,r})=\tilde{B}^\bullet_{r,j}(\tilde{B}^\bullet_{j,r}\tilde{B}^\bullet_{r,j})$}
\end{align}

\noindent
Then using $\tilde{B}^\bullet_{r,j}\tilde{B}^\bullet_{i,r}=\tilde{B}^\bullet_{j,r}\tilde{B}^\bullet_{r,i}$ from \eqref{eq:(***)} on the left of \eqref{eq:(v*)}, we have
$$\tilde{B}^\bullet_{j,r}\tilde{B}^\bullet_{r,i}\tilde{B}^\bullet_{r,i}=\tilde{B}^\bullet_{r,j}\tilde{B}^\bullet_{j,r}\tilde{B}^\bullet_{r,j}.$$
Cancelling $\tilde{B}^\bullet_{j,r}$, we have $$\tilde{B}^\bullet_{r,i}\tilde{B}^\bullet_{r,i}=\tilde{B}^\bullet_{r,j}\tilde{B}^\bullet_{r,j},$$
implying $\tilde{B}^\bullet_{r,i}=\tilde{B}^\bullet_{r,j}$ mod 2, which contradicts \eqref{eq:(+)}.

To complete the proof of ($*$), we need to show that it holds 
for $B_0$ such that $\Gamma(B_0)$ is a minimal infinite type diagram which is mutation-equivalent an extended Dynkin diagram of type $\tilde{C}_{n-1}$ (recall that $\Gamma(B_0)$ has $n$ vertices by our assumption; according to the convention in Figure~\ref{fig:extended-dynkin-diagrams} the diagram $\tilde{C}_{n-1}$ has $n$ vertices). Then it follows from the classification of minimal infinite type diagrams in \cite[Table 1]{S3} that $\Gamma(B_0)$ itself is an orientation of the extended Dynkin diagram $\tilde{C_{n-1}}$ (Figure~\ref{fig:extended-dynkin-diagrams}). Let us assume that the vertices of $\Gamma=\Gamma(B_0)$ are labeled $1,..,n$ in linear order (so the edges whose weights are equal to $2$ are $\{1,2\}$ and $\{n-1,n\}$). Let us note that in the diagram of $\tilde{\Gamma}=\Gamma(\tilde{B_0})$, each vertex $n+1\leq i \leq 2n$ is connected to exactly one vertex in $\Gamma(B_0)$. In particular, the vertex $2n$ is connected to exacly one vertex $n$ in $\Gamma(B_0)$. We denote by $B_0^n$ the $(n+1)\times n$ matrix obtained from $\tilde{B_0}$ by removing the rows indexed by $n+1,...,2n-1$; so its diagram $\Gamma(B_0^n)$ is the subdiagram consisting of $\Gamma$ and the vertex $2n$. 

We will show that ${B_0^n}$, so $\tilde{B_0}$, is of infinite mutation type using an inductive argument. For this let us note that for any $1\leq k \leq n-1$, the diagram $\mu_k(\Gamma(B_0^n))$ has a subdiagram $X$ which contains the edge $\{n-1,n\}$ such that $X$ is of type $\tilde{C}$ with the vertex $2n$ being connected to exactly one vertex, which is $n$, in $X$. Furthermore, if $k$ is not a source nor a sink in $\Gamma(B_0^n)$, then $X$ is a proper subdiagram. Applying mutations to the $\tilde{C}$ type diagrams obtained this way, we see that $B_0^n$ is mutation-equivalent to $B^n$ by a sequence of mutations at vertices from $1,2,...,n-1$ such that $\Gamma(B^n)$ contains a subdiagram $X'=\{i,n-1,n\}$, $i<n-1$, which is of type $\tilde {C}_2$ with the vertex $2n$ being connected to exactly one vertex, which is $n$, in $X'$. Then, applying a mutation at the vertex $i$ if necessary, we may assume that the vertex $n-1$ is not a source nor sink in $X'$. Then, in the diagram $\mu_{n-1}(B^n)$, the subdiagram $\{i,n\}$ has weight four and the vertex $2n$ is connected to exactly one vertex there. Then $\mu_{n-1}(B^n)$, so $\tilde{B_0}$, is of infinite mutation type by Proposition~\ref{pr:4-infinite}. This completes the proof of the theorem.

\end{document}